\documentclass{article}

\usepackage[hyphens]{url}
\usepackage{graphicx}
\usepackage{amsmath}
\usepackage{amsthm}
\usepackage{amsfonts}
\usepackage{amssymb}

\theoremstyle{theorem}
\newtheorem{theorem}{Theorem}
\newtheorem{lemma}[theorem]{Lemma}

\begin{document}

\begin{center}
	\vspace*{0.1cm}
	{\huge Taffy, Trees, and Tangles\par}
	\vspace{4ex}
	{\large Neil J. Calkin, Eliza Gallagher, and Ben Gobler\par\vspace{2\baselineskip}
	}
\end{center}

\begin{abstract}
We study the relationship between three combinatorial objects---a taffy pulling machine, the Calkin-Wilf tree of all fractions, and Conway's rational tangles. After introducing these objects, we develop a taffy analogue for Conway's characterization of rational tangles and we give a direct geometric connection between rational tangles and taffy pulls.
\end{abstract}

\section{Introduction}

Having multiple ways to represent the same mathematical concept often has advantages. Take for instance the problem of dissecting a square into smaller unequal squares, which gained notoriety in the early 20th century \cite{history}. Initial searches did not yield any solutions, leading many to question whether such a dissection was possible. When then-college students Brooks, Smith, Stone, and Tutte became interested in the problem, they realized that square dissections correspond to certain electrical networks, allowing them to apply established theory from electrical engineering to discover new construction methods \cite{brooks}. Insight from this second representation of the problem led to a breakthrough in finding the first ``squared square" consisting of 69 tiles \cite{stewart}.

Each representation of a concept may offer its own advantages, such as notational ease or conceptual simplicity. If nothing else, the joy of connecting multiple representations will make a mathematician smile. We hope to foster that experience by connecting three combinatorial objects from the literature: taffy pulling machines, the Calkin-Wilf tree of all fractions, and Conway's rational tangles.

\section{Stretch and Fold}

Taffy is a candy made from a heated sugar mixture. As it cools, the taffy is \textit{pulled}, a kneading process that incorporates air to make the taffy soft and chewy. Taffy pulling can be done by hand or with a machine; there are numerous designs of taffy pulling machines that are used commercially and have various mathematical properties \cite{thiffeault}.

\begin{figure}[!hbt]
\centering
\includegraphics[width=\textwidth]{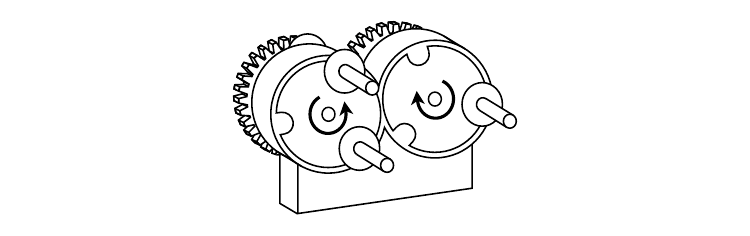}
\caption{A three-peg taffy pulling machine.}
\label{taffy_machine}
\end{figure}

\begin{figure}[!hbt]
\centering
\includegraphics[width=\textwidth]{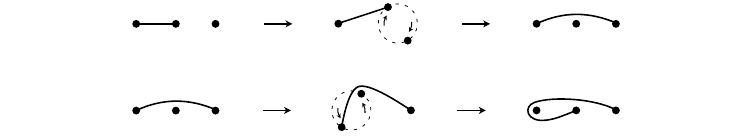}
\caption{Visualizing a right turn followed by a left turn.}
\label{turn_demos}
\end{figure}

In this paper, we consider the design patented by Nitz in 1918 (Figure \ref{taffy_machine}). A Nitz taffy machine has three large pegs that rotate around each other in figure eight motions. A confectioner loads the machine by attaching a wad of taffy to the leftmost two pegs. When the machine runs, the motion of the pegs can be simplified into two actions: a \textit{right turn}, in which the rightmost two pegs trade places by rotating clockwise, and a \textit{left turn}, in which the leftmost two pegs trade places by rotating counterclockwise (Figure \ref{turn_demos}). In each case, the middle peg rotates upward and outward.

\begin{figure}[!hbt]
\centering
\includegraphics[width=\textwidth]{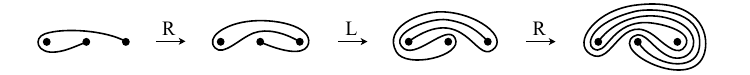}
\caption{Depicting a sequence of taffy pulls made by the Nitz machine.}
\label{alternating_turns}
\end{figure}

We define a \textit{taffy pull} to be any finite sequence of left and right turns in any order. For example, we write $t = RRL = R^2L$ to denote the taffy pull made by two right turns followed by a left turn. We include the sequence of zero turns and call it the \textit{initial taffy pull}, denoted $t^0$. The Nitz machine is designed to make alternating left and right turns. Continuing from Figure \ref{turn_demos}, Figure \ref{alternating_turns} illustrates how the taffy folds as the machine makes three more turns.

Our \textit{taffy diagrams} are cross-sectional views of the taffy machine. To balance the representational accuracy with visual clarity, we abide by a few conventions. Strands of taffy that touch each other will stick together in practice, but for better visibility we leave small gaps between nearby strands. To preserve the appearance of the original taffy pull, we avoid unnecessary slack, and we do not allow the strands of taffy to cross each other or pass through the pegs.

The aim of pulling taffy is to trap air bubbles between the many folds. The amount of trapped air is roughly proportional to the number of \textit{layers} of taffy, where a layer refers to a horizontal segment of taffy between adjacent pegs. An effective taffy machine produces increasingly many layers of taffy with each turn. To measure the effectiveness of the Nitz machine, we calculate the total number of layers of taffy after $n$ turns and determine the growth rate of this quantity with respect to $n$ (Thiffeault similarly defines the efficiency of a taffy machine in \cite{thiffeault}).

\begin{figure}[!hbt]
\centering
\includegraphics[width=\textwidth]{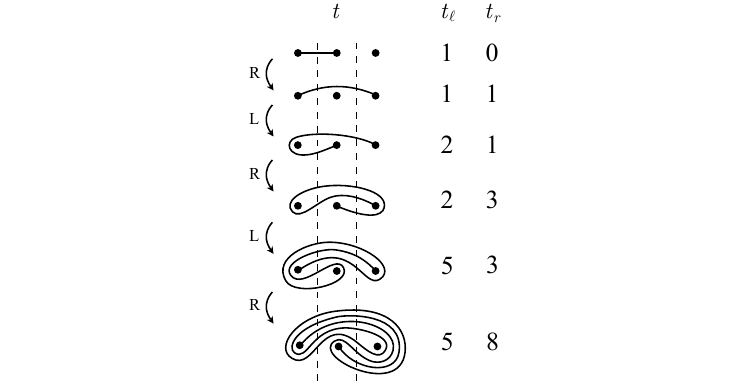}
\caption{Counting the numbers of left and right layers.}
\label{layer_counts}
\end{figure}

For the purpose of counting, it helps to distinguish the layers of taffy appearing between the leftmost two pegs as \textit{left layers}, and those appearing between the rightmost two pegs as \textit{right layers}. For a taffy pull $t$ we write $t_\ell$ to denote the number of left layers and $t_r$ to denote the number of right layers. These numbers are obtained from the diagram of $t$ by drawing vertical lines between adjacent pegs. Then $t_\ell$ counts the number of times that the taffy crosses the line between the leftmost pegs, and $t_r$ counts the number of times that the taffy crosses the line between the rightmost pegs. Counting the left and right layers of taffy pulls made by the Nitz machine reveals a familiar pattern (Figure \ref{layer_counts}). Indeed, the layer counts are all Fibonacci numbers: after $n$ turns, the numbers of left and right layers are $F_{n+1}$ and $F_n$ respectively if $n$ is even, or $F_n$ and $F_{n+1}$ if $n$ is odd. To see why this pattern continues, we will determine in general how each type of turn affects the layer counts.

\begin{figure}[!hbt]
\centering
\includegraphics[width=\textwidth]{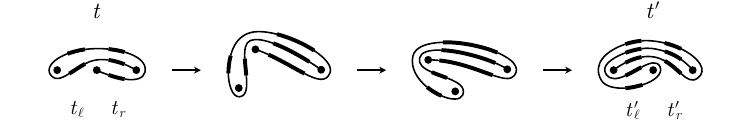}
\caption{A left turn ``adds" the right layers to the left layers.}
\label{adding_layers}
\end{figure}

Let $t$ be a taffy pull with $t_\ell$ left layers and $t_r$ right layers. When $t$ is extended by a left turn, the right layers of $t$ get stretched to the left, and the left layers of $t$ rotate underneath. This is demonstrated in Figure \ref{adding_layers}, where layers are shown with a thickened middle section for clarity in visualization. When the turn is complete, the resulting taffy pull $t' = tL$ has layer counts $t'_\ell = t_\ell + t_r$ and $t'_r = t_r$.

\newpage

\noindent Similarly, when $t$ is extended by a right turn, the resulting taffy pull $t' = tR$ has layer counts $t'_\ell = t_\ell$ and $t'_r = t_\ell + t_r$.

We apply this result to count the layers in taffy pulls made by the Nitz machine. The initial taffy pull has $1 = F_1$ left layers and $0 = F_0$ right layers. The machine makes a right turn, and the resulting taffy pull still has $F_1$ left layers but now has $F_2 = F_1 + F_0$ right layers. Then a left turn is made, and the resulting taffy pull has $F_3 = F_2 + F_1$ left layers and $F_2$ right layers. Since $F_n = F_{n-1} + F_{n-2}$, it follows that after $n$ alternating turns, the numbers of left and right layers are $F_{n+1}$ and $F_n$ respectively if $n$ is even, or $F_n$ and $F_{n+1}$ if $n$ is odd.

Now we measure the effectiveness of the Nitz machine. After $n$ turns, the taffy pull $t$ has a total number of layers $t_\ell + t_r = F_n + F_{n+1} = F_{n+2}$, again a Fibonacci number (Margalit observes this pattern in \cite{margalit}). Hence, the $n$th turn increases the total number of layers by a factor of $F_{n+2}/F_{n+1}$, which approaches the golden ratio $\varphi = (1+\sqrt{5})/2 \approx 1.618$ as $n$ grows. Accordingly, the total number of layers after $n$ turns is asymptotic to $\varphi^n/\sqrt{5}$. Exponential growth of layer counts is common among conventional taffy machine designs, but the growth factor varies between them. Is $F_{n+2}$ total layers after $n$ turns optimal for a three-peg taffy machine? To answer this question, we consider other turning sequences combinatorially.

\section{Beyond Nitz}

The Nitz machine's alternating turns produce taffy pulls with Fibonacci numbers of layers. What about other sequences of turns? We could imagine a taffy machine that can make any sequence of left and right turns upon request, and so can produce any taffy pull. It is reasonable to ask which pairs of integers can be the left and right layer counts of some taffy pull.

\begin{figure}[!hbt]
\centering
\includegraphics[width=\textwidth]{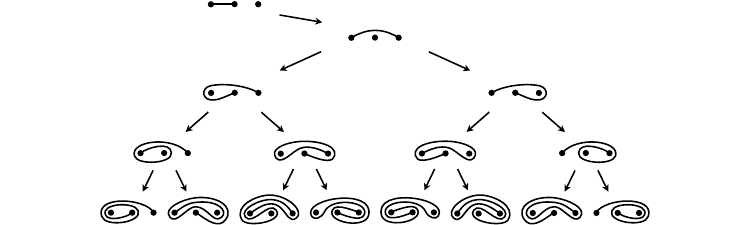}
\caption{The tree of taffy pulls, illustrated with taffy diagrams.}
\label{taffy_tree}
\end{figure}

Taffy pulls can be arranged in a binary tree where each node has two children, called the \textit{left child} and the \textit{right child} (Figure \ref{taffy_tree}). The left and right children of a node labeled by the taffy pull $t$ are the extensions of $t$ by a single left or right turn respectively. Extending the initial taffy pull $t^0$ by a left turn does not change its diagram, so we omit the left branch of the tree for now. The taffy diagrams in Figure \ref{taffy_tree} are distinct, suggesting that each taffy pull beginning with a right turn has a distinct diagram---we will return to this shortly.

\newpage
We denote the downward path in the tree to the taffy pull $t = R^{n_1}L^{n_2}\cdots R^{n_k}$ using the same notation, $R^{n_1}L^{n_2}\cdots R^{n_k}$. And for each taffy pull $t$, we write the ordered pair $(t_r,t_\ell)$ to denote the numbers of right and left layers of $t$; the seemingly reversed order is intended to make the following connection more intuitive. We determined in the previous section that the left child of a node labeled with the taffy pull $t$ has layer counts $(t_r,t_\ell+t_r)$, and the right child has layer counts $(t_\ell+t_r,t_\ell)$. We hope that some readers may recognize this pair of generating actions from what has become known as the \textit{Calkin-Wilf tree}: a binary tree of fractions with root $1/1$ in which the left child of $a/b$ is $a/(a+b)$ and the right child is $(a+b)/b$ \cite{calkin}. We call these generating actions the \textit{left rule} and \textit{right rule} respectively.

\begin{figure}[!hbt]
\centering
\includegraphics[width=\textwidth]{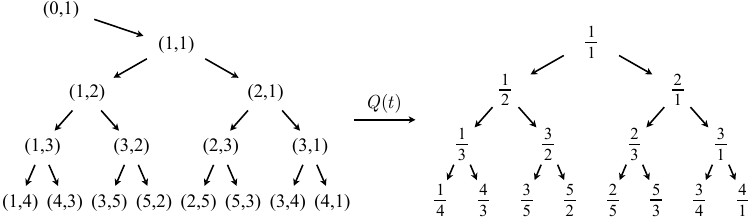}
\caption{Layer counts $(t_r,t_\ell)$ and the Calkin-Wilf tree.}
\label{cw_tree}
\end{figure}

To make the correspondence explicit, we assign each taffy pull $t$ with layer counts $(t_r,t_\ell)$ the \textit{taffy number} $Q(t) = t_r / t_\ell$. Then $t$ and $Q(t)$ are located at the same nodes in the tree of taffy pulls and the Calkin-Wilf tree respectively, as seen in Figure \ref{cw_tree}.

The Calkin-Wilf tree has two characteristics that are particularly useful in this context: every fraction in the tree is in lowest terms, and every positive rational number appears exactly once in the tree. Consequently, every taffy pull $t$ has layer counts $t_\ell$ and $t_r$ that are coprime, and for any pair of coprime positive integers $(a,b)$, there is exactly one taffy pull $t$ beginning with a right turn such that $(t_r,t_\ell) = (a,b)$. Hence, taffy pulls beginning with a right turn have distinct taffy diagrams, as suggested by Figure \ref{taffy_tree}. In the next section we will see that the sequence of left and right steps to $a/b$ in the Calkin-Wilf tree can be obtained by applying the Euclidean algorithm to $a$ and $b$, and therefore the taffy pull with a given taffy number $a/b$ can be located in the tree of taffy pulls in time logarithmic in the size of $a$ and $b$.

Returning to the discussion of effectiveness, we can compare the Nitz machine's alternating turns to any other sequence of left and right turns. For a positive integer $n$, we are interested in determining which sequence of $n$ turns produces a taffy pull $t$ with the most total layers $t_\ell + t_r$, since this number is roughly proportional to the amount of trapped air in the taffy. The largest total corresponds to the largest numerator-denominator sum among all fractions in row $n$ of the Calkin-Wilf tree, which, in turn, is the largest numerator (or denominator) in row $n+1$. The largest numerator in row $n$ of the Calkin-Wilf tree is the Fibonacci number $F_{n+1}$ \cite{wilf}. Hence, the largest total number of layers in a taffy pull made by $n$ turns is $F_{n+2}$, showing that the alternating turns performed by the Nitz machine are in fact the most effective.

\section{Tracing History}

While taffy diagrams are a visually appealing way to represent taffy pulls, taffy numbers are a more useful representation for answering certain questions. In the previous section we found that all taffy pulls beginning with a right turn have distinct diagrams. Given one of these diagrams, we may wish to recover the sequence of left and right turns that generated it.

We see from Figure \ref{taffy_tree} that the diagrams of taffy pulls below $t=R$ are lopsided, having more mass to the left or to the right of the center peg depending on the final turn. With some effort, the diagram can be redrawn with the final turn undone. This process can be repeated to identify the penultimate turn, and so on, eventually obtaining the sequence of turns that generated the diagram. Each iteration makes a single step up the tree of taffy pulls, and doing so is rather laborious.

Changing representations allows us to work directly with taffy numbers, which turns out to be far more convenient. In the Calkin-Wilf tree, $a/b$ is a left child if $a < b$, and a right child if $a > b$. The parent of $a/b$ is found by subtracting the smaller of $a$ and $b$ from the larger. This procedure for stepping up the Calkin-Wilf tree one level at a time is essentially the same as the original algorithm in Proposition 2 of \textit{Euclid's Elements Book VII} for finding the greatest common divisor of $a$ and $b$, sometimes called the \textit{slow Euclidean algorithm}, which proceeds by subtraction rather than division \cite{Euclid}.

\begin{figure}[!hbt]
\centering
\includegraphics[width=\textwidth]{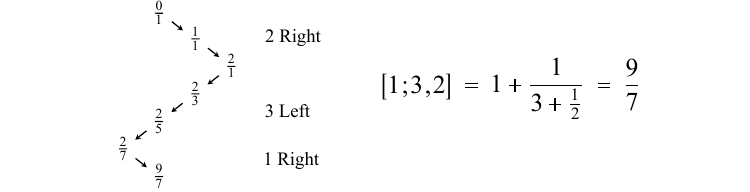}
\caption{The path from $0/1$ to $9/7$ corresponds to the continued fraction for $9/7$.}
\label{path_cf}
\end{figure}

Writing taffy pulls in terms of consecutive runs of the same type of turn provides a connection to the better known \textit{fast Euclidean algorithm}, which replaces the larger of $a$ and $b$ with the remainder upon division by the smaller number. The integer quotient counts the maximal number of consecutive subtractions of one number from the other, so each iteration of the algorithm slides up the Calkin-Wilf tree along a maximal run of steps in the same direction. Further, the quotients obtained from successive iterations of the fast Euclidean algorithm are the coefficients of a continued fraction for $a/b$. Hence, the path from $0/1$ to $a/b$ in the Calkin-Wilf tree corresponds to the continued fraction for $a/b$ whose coefficients are the numbers of consecutive steps in each direction along the path \cite{gobler}. Explicitly, the path $R^{n_1}L^{n_2}\cdots R^{n_k}$ ends at the rational number with continued fraction representation $[n_k;n_{k-1},\dots,n_1]$. Figure \ref{path_cf} illustrates the path from $0/1$ to $9/7$ corresponding to the continued fraction for $9/7$.

The notion of backward steps in the Calkin-Wilf tree turns out to be highly profitable, as we will explore in the next two sections.

\section{Working backward}

\begin{figure}[!hbt]
\centering
\includegraphics[width=\textwidth]{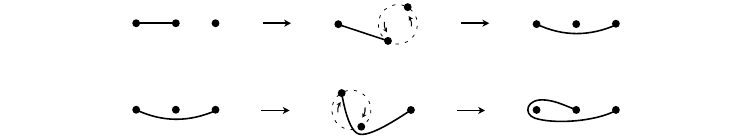}
\caption{Visualizing a reverse right turn followed by a reverse left turn.}
\label{reverse_turns}
\end{figure}

Imagine that a switch on the back of the taffy machine reverses the direction of the motor---now the leftmost two pegs rotate clockwise, and the rightmost two pegs rotate counterclockwise. In each case, the middle peg moves downward and outward, rather than upward and outward. Figure \ref{reverse_turns} illustrates a \textit{reverse right turn} followed by a \textit{reverse left turn}.

Since each reverse turn undoes its corresponding forward turn, we denote a reverse left turn and reverse right turn by $L^{-1}$ and $R^{-1}$ respectively. As a counterpart to the original taffy pulls (now more appropriately called \textit{forward taffy pulls}), we define a \textit{reverse taffy pull} to be any finite sequence of reverse left and reverse right turns.

\begin{figure}[!hbt]
\centering
\includegraphics[width=\textwidth]{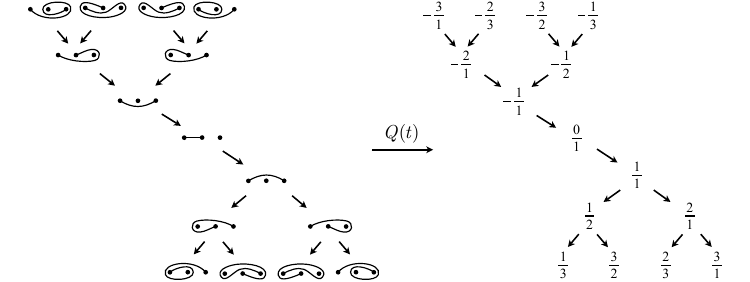}
\caption{The double tree of fractions assigns a negative taffy number to each reverse taffy pull.}
\label{double_trees}
\end{figure}

Naturally, we wish to assign a taffy number to each reverse taffy pull. This can be done following the observation that the left and right rules generating the Calkin-Wilf tree are invertible: the parent of a left child $a/b$ is $a/(b-a)$, and the parent of a right child $a/b$ is $(a-b)/b$. We call these operations the \textit{backward left rule} and the \textit{backward right rule} respectively. Applying the backward right rule to the root of the tree $1/1$ yields $0/1$, which is the taffy number of the initial taffy pull. Applying the backward right rule again yields $-1/1$. From there, additional backward left and right steps produce a negative analog of the Calkin-Wilf tree. We call this extension of the Calkin-Wilf tree the \textit{double tree} of fractions (Figure \ref{double_trees}). Accordingly we assign the reverse taffy pull $t = R^{-n_1}L^{-n_2}\cdots R^{-n_k}$ the taffy number $Q(t)$ given by the fraction at the end of the same path in the double tree (where we are careful to note that a backward right step is made by moving up and to the left, while a backward left step is made by moving up and to the right). It can be seen that the reverse taffy pull $t$ with layer counts $(t_r,t_\ell)$ has taffy number $Q(t) = -t_r/t_\ell$.

Each rational number appears exactly once in the double tree of fractions. Hence, the double tree of taffy pulls has a representative for each taffy number. We call a taffy pull in the double tree the \textit{canonical} taffy pull for its taffy number; these are the forward and reverse taffy pulls beginning with $R$ or $R^{-1}$, as well as $t^0$.

The diagrams of canonical taffy pulls as seen in Figure \ref{double_trees} are related by rotational symmetry. If the diagrams of $t$ and $t'$ are $180^\circ$ rotations of each other, then $Q(t)$ and $Q(t')$ are negative reciprocals of each other. For instance, the canonical taffy pulls $t = RLR$ and $t' = R^{-2}L^{-1}$ have diagrams that are 180$^\circ$ rotations of each other, and their taffy numbers are $Q(t) = 3/2$ and $Q(t') = -2/3$. We will use this fact to prove Lemma \ref{single_extension} in the next section, but first we must introduce a notational convenience.

No taffy pull has a diagram that is the 180$^\circ$ rotation of the diagram of the initial taffy pull $t^0$. In lieu of one, we use $t^{\infty}$ as a symbolic placeholder for a taffy pull whose diagram is the 180$^\circ$ rotation of the diagram of $t^0$, in which the taffy connects the rightmost two pegs instead of the leftmost two pegs. The taffy number of $t^\infty$ is then $Q(t^\infty) = -1/Q(t^0) = -1/0 = 1/0$. We let $t^\infty$ be the canonical taffy pull for the taffy number $1/0$. Now, the 180$^\circ$ rotated diagram of each canonical taffy pull is the diagram of another canonical taffy pull. This is the infrastructure that we need for the next section, where we will generalize what is considered to be a taffy pull.

\section{Forward and Backward} 

As we have seen, the freely turning taffy machine can make a taffy pull with any rational number as its taffy number by operating either normally or in reverse. What about a machine that can make both forward and reverse turns in any order?

\begin{figure}[!hbt]
\centering
\includegraphics[width=\textwidth]{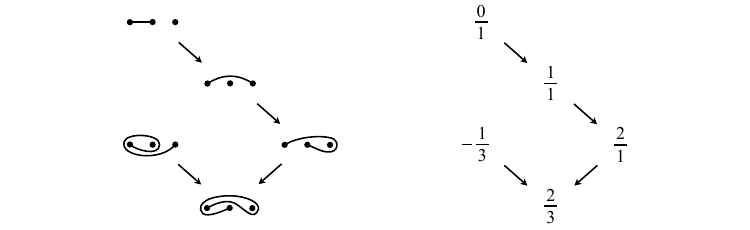}
\caption{The generalized taffy pull $t = R^2LR^{-1}$ has taffy number $-1/3$.}
\label{taffy_example}
\end{figure}

We define a \textit{generalized taffy pull} to be any finite sequence of left, right, reverse left, or reverse right turns. For example, $t = R^2LR^{-1}$ is a generalized taffy pull that is neither a forward taffy pull nor a reverse taffy pull. We wish to assign a taffy number to each generalized taffy pull, and we already have a system in place: starting with $0/1$, apply the left, right, backward left, and backward right rules corresponding to each turn. Following this convention, $t = R^2LR^{-1}$ has taffy number $-1/3$ (Figure \ref{taffy_example}).

Since the pairs of turns $LL^{-1}$,  $RR^{-1}$, $L^{-1}L$, and $R^{-1}R$ have no net effect on the taffy machine or the taffy number, we adopt the convention that such pairs will be collapsed, so that, for example, $t = R^{-1}R^2L$ is written more simply as $t = RL$. Then every taffy pull can be written in the standard form $t = R^{n_1}L^{n_2}\cdots R^{n_k}$ where each $n_i$ is a nonzero integer for $i=2,\dots,k-1$ ($n_1$ or $n_k$ may be zero if $t$ begins or ends with a left turn); $n_i > 0$ indicates forward turns and $n_i < 0$ indicates reverse turns.

\begin{figure}[!hbt]
\centering
\includegraphics[width=\textwidth]{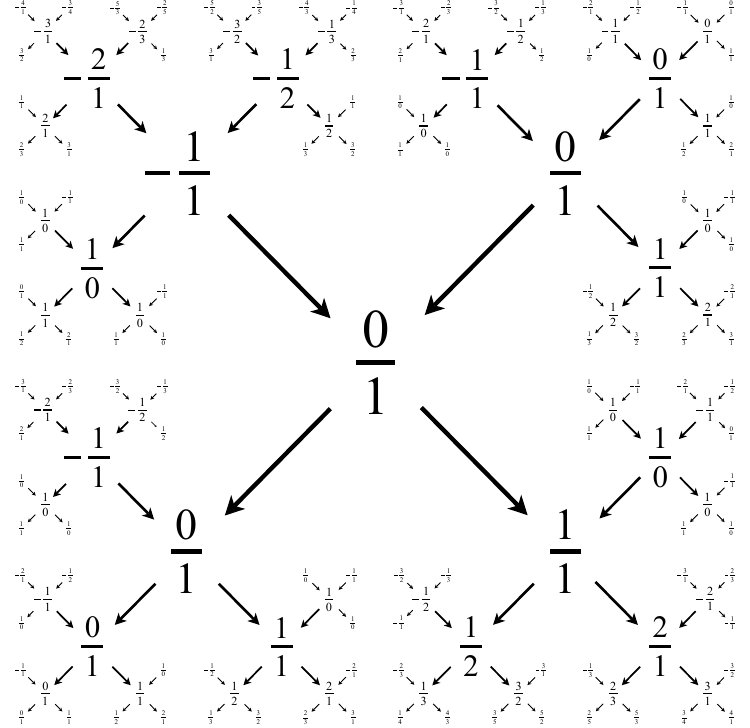}
\caption{The four-way tree of fractions.}
\label{four_way}
\end{figure}

Allowing forward and reverse turns in any order gives rise to an extension of the Calkin-Wilf tree in which forward and backward left and right steps are made from every node. The result is what we call the \textit{four-way tree} of fractions (Figure \ref{four_way}, \cite{gobler}), which contains the taffy number of every generalized taffy pull. We are again careful to note that backward right steps are made by moving up and to the left, while backward left steps are made by moving up and to the right.

Every fraction, including $1/0$, appears infinitely often in the four-way tree. In order for $t^\infty$ to be well-defined, we select $t^\infty = RL^{-1}$ to be the canonical taffy pull for the taffy number $1/0$, noting that the diagram of $RL^{-1}$ is indeed the $180^\circ$ rotation of the diagram of $t^0$.

When two generalized taffy pulls have the same taffy number, how do their diagrams compare? For example, the generalized taffy pull $t = R^2LR^{-1}$ shown in Figure \ref{taffy_example} has taffy number $-1/3$, and the canonical taffy pull for the taffy number $-1/3$ is $t^* = R^{-1}L^{-2}$, whose diagram appears in the upper-right corner of the double tree in Figure \ref{double_trees}. Remarkably, the diagrams of $t$ and $t^*$ are identical. We call two generalized taffy pulls \textit{equivalent} if they have the same diagram. We will show that two generalized taffy pulls are equivalent if and only if they have the same taffy number.

\begin{figure}[!hbt]
\centering
\includegraphics[width=\textwidth]{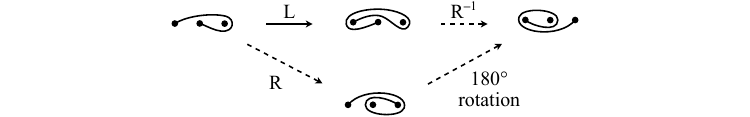}
\caption{The pair of turns $LR^{-1}$ has the same effect as a right turn followed by a 180$^\circ$ rotation.}
\label{critical_pair}
\end{figure}

\begin{lemma} \label{equivalences} Let $t$ be a taffy pull. Then the diagram of $tLR^{-1}$ is the $180^\circ$ rotation of the diagram of $tR$.
\end{lemma}
\begin{proof}
The composition of turns $LR^{-1}$ rotates the leftmost peg underneath both pegs to its right. This same motion of the pegs is achieved by a right turn followed by a $180^\circ$ rotation (Figure \ref{critical_pair}). Hence, the diagram of $tLR^{-1}$ is the $180^\circ$ rotation of the diagram of $tR$.
\end{proof}

Analogous statements hold by a similar analysis: the diagrams of $tRL^{-1}$, $tL^{-1}R$, and $tR^{-1}L$ are the $180^{\circ}$ rotations of the diagrams of $tL$, $tR^{-1}$, and $tL^{-1}$ respectively.

\begin{lemma} \label{single_extension}
Let $t$ be a canonical taffy pull, and let $s$ be a single turn. Then the taffy pull $ts$ is equivalent to the canonical taffy pull with the same taffy number as $ts$.
\end{lemma}
\begin{proof}
If $t = t^0$ or $t^\infty$, it is quickly verified that extending $t$ by any single turn $s$ results in a taffy pull $ts$ that has the same diagram as the canonical taffy pull $t^*$ with the same taffy number, either $t^* = t^0$, $t^\infty$, $R$, or $R^{-1}$. Next, let $t$ be any other canonical taffy pull. We consider the case where $t$ is a forward taffy pull. Then $tL$ and $tR$ are canonical taffy pulls, for which the result is immediate. To address $tL^{-1}$ and $tR^{-1}$, we consider further the case where $t$ ends with a left turn and write $t$ as $cL$, where $c$ is a canonical forward taffy pull. Then $tL^{-1} = cLL^{-1} = c$ is canonical, for which the result is immediate. By Lemma \ref{equivalences}, the diagram of $tR^{-1} = cLR^{-1}$ is the $180^\circ$ rotation of the diagram of $cR$. As $cR$ is a canonical taffy pull, there is a canonical taffy pull $t^*$ whose diagram is the $180^\circ$ rotation of the diagram of $cR$, and hence is equivalent to $tR^{-1}$. It remains to show that $tR^{-1}$ and $t^*$ have the same taffy number. Let $Q(c) = a/b$. Then $Q(tR^{-1}) = Q(cLR^{-1}) = -b/(a + b)$ by applying the left rule to $Q(c)$ followed by the right rule. Since pairs of canonical taffy pulls with 180$^\circ$ rotated diagrams have negative reciprocal taffy numbers, $Q(t^*) = -1/Q(cR) = -1/((a + b)/b) = -b / (a + b)$. Thus, $Q(tR^{-1}) = Q(t^*)$, showing $tR^{-1}$ is equivalent to the canonical taffy pull $t^*$ with the same taffy number. Similar arguments hold when instead $t$ is a forward taffy pull that ends with $R$, or when $t$ is a reverse taffy pull that ends with $L^{-1}$ or $R^{-1}$.
\end{proof}

We extend the result of Lemma \ref{single_extension} by induction to all generalized taffy pulls.

\begin{lemma} \label{equivalent_canonical}
Let $t$ be a generalized taffy pull. Then $t$ is equivalent to the canonical taffy pull with the same taffy number.
\end{lemma}
\begin{proof}
We proceed by induction on $n(t)$, the number of turns comprising $t$. If $n(t) = 0$, then $t = t^0$ is canonical, so the result is immediate. If $n(t)\geq1$, write $t = t's$ where $s$ is a single turn. By the inductive hypothesis, $t'$ is equivalent to the canonical taffy pull $(t')^*$ with the same taffy number. Extending both by the turn $s$, $t's$ is equivalent to $(t')^*s$ and has the same taffy number. By Lemma \ref{single_extension}, $(t')^*s$ is equivalent to the canonical taffy pull $t^*$ with the same taffy number. Hence, $t = t's$ is equivalent to $t^*$ and has the same taffy number.
\end{proof}

\begin{theorem} \label{equivalent_generalized}
Two generalized taffy pulls are equivalent if and only if they have the same taffy number.
\end{theorem}
\begin{proof}
Let $t_1,t_2$ be generalized taffy pulls. By Lemma \ref{equivalent_canonical}, there are canonical taffy pulls $t_1^*,t_2^*$ equivalent to $t_1,t_2$ with the same taffy numbers respectively. Then $t_1,t_2$ are equivalent $\iff$ $t_1^*,t_2^*$ are equivalent $\iff$ $t_1^*,t_2^*$ have the same taffy number $\iff$ $t_1,t_2$ have the same taffy number.
\end{proof}

Taffy pulls are not alone in having a classification like this. In particular, the statement of Theorem \ref{equivalent_generalized} has the same flavor as Conway's theorem on the classification of rational tangles. In the next section, we give an overview of the theory of rational tangles that is similar to our development of taffy pulls.

\section{Getting Tangled}

\begin{figure}[!hbt]
\centering
\includegraphics[width=\textwidth]{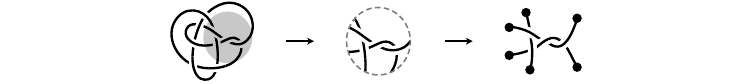}
\caption{A three-arc tangle found within a knot diagram.}
\label{knot_tangle}
\end{figure}

The notion of a tangle in knot theory is used to describe the interactions between strands in a knot diagram. By wrapping a sphere around a portion of a knot, the interior space isolates a number of arcs with their ends on the surface of the sphere (Figure \ref{knot_tangle}). Such a collection of arcs is called a \textit{tangle}.

In 1970, J. H. Conway introduced a family of two-arc tangles called \textit{rational tangles} to provide a convenient notation for the enumeration of knots and links \cite{conway}. More recently, rational tangles have become a valuable representation in molecular biology for studying DNA topology \cite{goldman}. In this section we consider rational tangles as combinatorial objects with a connection to the Calkin-Wilf tree.

\begin{figure}[!hbt]
\centering
\includegraphics[width=\textwidth]{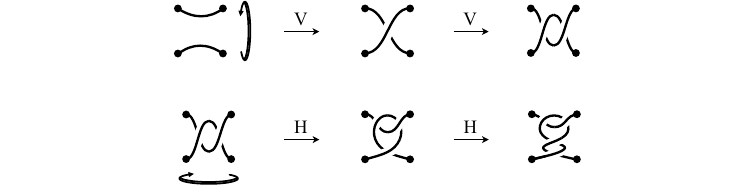}
\caption{Two vertical twists followed by two horizontal twists.}
\label{tangle_twists}
\end{figure}

A rational tangle begins as a pair of horizontal arcs with their ends positioned on the four corners of a square. The ends may be twisted around each other in one of two ways: a \textit{vertical twist}, made by lifting the lower-right end in front of the upper-right end as they trade positions, and a \textit{horizontal twist}, made by lifting the lower-right end in front of the lower-left end as they trade positions (Figure \ref{tangle_twists}).

\begin{figure}[!hbt]
\centering
\includegraphics[width=\textwidth]{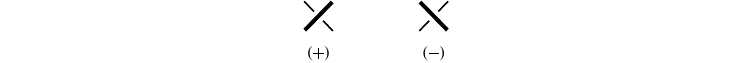}
\vspace{-0.6cm}
\caption{A positive crossing and a negative crossing.}
\label{crossing_signs}
\end{figure}

Vertical and horizontal twists can be undone; we call the reverse motions a \textit{vertical untwist} and \textit{horizontal untwist} respectively. In a tangle diagram, twists and untwists are distinguished by the crossings they form: each twist, whether vertical or horizontal, forms a \textit{positive crossing}, and each untwist forms a \textit{negative crossing} (Figure \ref{crossing_signs}) \cite{adams}.

\begin{figure}[!hbt]
\centering
\includegraphics[width=\textwidth]{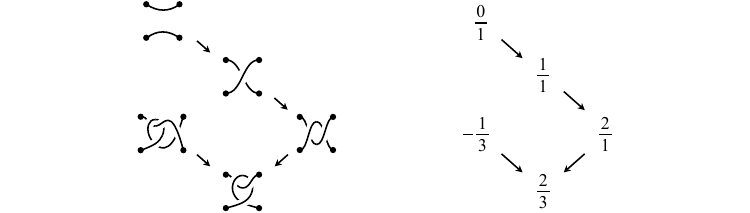}
\caption{The rational tangle $T = V^2HV^{-1}$ has tangle number $-1/3$.}
\label{tangle_example}
\end{figure}

We define a \textit{rational tangle} to be any finite sequence of vertical or horizontal twists or untwists. By now it should come as no surprise that the four types of twists can be associated with the directions of movement in the four-way tree in Figure \ref{four_way}. To each tangle $T = V^{n_1}H^{n_2}\cdots V^{n_k}$ we assign a \textit{tangle number} $Q(T)$ given by the fraction at the end of the path $R^{n_1}L^{n_2}\cdots R^{n_k}$ in the four-way tree. From the discussion in Section 4, it follows that $Q(T) = [n_k;n_{k-1},\dots,n_1]$, the continued fraction whose coefficients are the numbers of consecutive steps in each direction along the path. For example, the rational tangle $T = V^2HV^{-1}$ has tangle number $Q(T) = [-1;1,2] = -1/3$ (Figure \ref{tangle_example}). This formulation of the tangle number in terms of a continued fraction is not only concise, but it is precisely how Conway assigned rational numbers to rational tangles in his original presentation \cite{conway}.

\begin{figure}[!hbt]
\centering
\includegraphics[width=\textwidth]{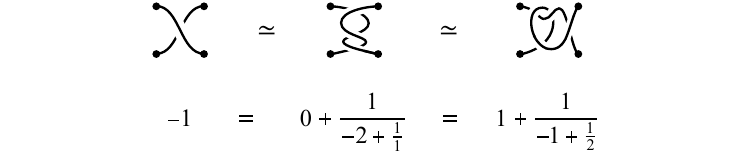}
\caption{Three equivalent rational tangles with tangle number $-1$.}
\label{equivalent_tangles}
\end{figure}

Conway's wonderful result is that two rational tangles have the same tangle number if and only if their diagrams are equivalent up to ambient isotopy (the diagrams of $T_1$ and $T_2$ are \textit{ambient isotopic} if the arcs of $T_1$ can be continuously deformed into those of $T_2$ without moving the endpoints and without the arcs passing through each other or leaving the square boundary formed by the endpoints \cite{conway}). Various proofs of this theorem are given in \cite{goldman} and \cite{kauffman}. To demonstrate, Figure \ref{equivalent_tangles} shows the diagrams of three rational tangles with the same tangle number, $-1$. By pulling the ends taut, all three tangles are seen to be equivalent.

\section{Twists and Turns}

We have surveyed two families of objects, taffy pulls and rational tangles, each with a natural correspondence to the four-way tree of fractions such that two objects have equivalent diagrams if and only if they are assigned the same rational number. In this sense, taffy pulls and rational tangles are two physical representations of the four-way tree of fractions. To change directly from one representation to the other, we showcase a geometric relationship between their diagrams.

\begin{figure}[!hbt]
\centering
\includegraphics[width=\textwidth]{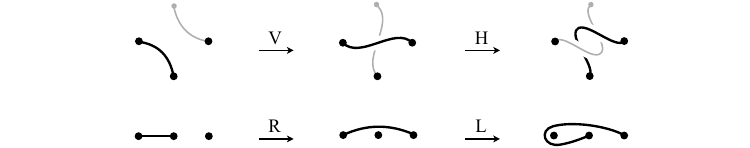}
\caption{Making twists and turns in tandem.}
\label{twists_turns}
\end{figure}

Taffy diagrams and rational tangle diagrams are structurally similar in that they both consist of points and arcs. Four of their features align once the tangle diagrams are rotated 45$^\circ$ clockwise: comparing the rotated initial rational tangle to the initial taffy pull, three endpoints of the tangle are roughly arranged like the pegs of the taffy machine, and one of the arcs is in the same position as the taffy (Figure \ref{twists_turns}). Making twists and turns in tandem, we explore whether those visual similarities persist.

\begin{figure}[!hbt]
\centering
\includegraphics[width=\textwidth]{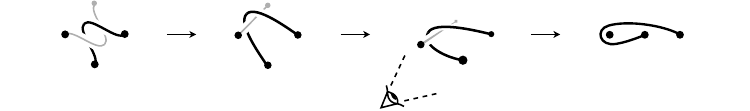}
\caption{Viewing the tangle diagram from a different perspective.}
\label{change_perspective}
\end{figure}

After a vertical twist and a right turn, the corresponding features in the tangle diagram still resemble those of the taffy machine, albeit somewhat abstractly. After a horizontal twist and a left turn, though, it is harder to see their likeness. Further, each tangle diagram has an additional arc and endpoint that do not immediately appear to have counterparts in taffy diagrams. But there is an epiphany about the roles of these features that comes about from a change in perspective. Figure \ref{change_perspective} transforms the tangle diagram by pulling taut the second arc and changing the viewing angle so that our eye sees directly along the straightened arc.

After the transformation, the second arc in the tangle diagram is positioned as one of the physical pegs of the taffy machine, an insight easy to miss when only viewing cross-sections of the taffy machine. From this observation, it is revealed that the taffy machine is truly a \textit{tangle} machine; the two arcs that get tangled together are the taffy and the third peg that is not attached to the taffy. How beautiful indeed!

\begin{figure}[!hbt]
\centering
\includegraphics[width=\textwidth]{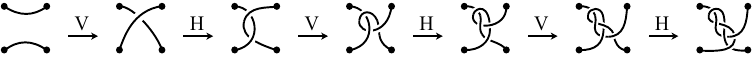}
\caption{Alternating vertical and horizontal twists produce a standard three-strand braid.}
\label{alternating_twists}
\end{figure}

While geometrically similar, taffy pulls and rational tangles offer different advantages in how they present information about their corresponding rational numbers. A taffy diagram readily shows the numerator and denominator of its taffy number via its layer counts, and the sign of the fraction is determined by whether the diagram is that of a forward or reverse taffy pull. These two pieces of information are not immediately visible in a tangle diagram, which instead displays the complete history of the twisting steps, something that taffy diagrams cannot do. Tangle diagrams are easier to draw since each new twist extends the previous diagram, unlike taffy diagrams which must be redrawn after each turn. And tangle diagrams tend to be more compact than taffy diagrams, whereas taffy diagrams become larger as the number of layers increases. As an example, we invite readers to draw the diagram of $t=RLRLRL$ (extending the alternating taffy pull in Figure \ref{alternating_turns} by a left turn), which has thirteen left layers and eight right layers. Meanwhile, the tangles produced by alternating vertical and horizontal twists form a visually compact three-strand braid (Figure \ref{alternating_twists}).

\section{Conclusion}

We have seen how taffy pulls, the four-way tree, and rational tangles are three representations of the same mathematical structure. Each representation has its own advantages: taffy pulls are conceptually simple, the four-way tree is mathematically convenient, and rational tangles are diagramatically compact. Do other combinatorial objects share the same structure as these three, and if so, what advantages might they offer?

\begin{center}
\textit{Oh, what a tangled web we weave}

\textit{when first we study taffy trees!}

\vspace{0.05in}\hspace{1.5in}\small -- With apologies to Sir Walter Scott
\end{center}

\bibliographystyle{amsplain}
\bibliography{References}

@article {thiffeault,
    AUTHOR = {Thiffeault, Jean-Luc},
     TITLE = {The mathematics of taffy pullers},
   JOURNAL = {Math. Intelligencer},
  FJOURNAL = {The Mathematical Intelligencer},
    VOLUME = {40},
      YEAR = {2018},
    NUMBER = {1},
     PAGES = {26--35},
      ISSN = {0343-6993,1866-7414},
   MRCLASS = {37E30 (70B15)},
  MRNUMBER = {3775163},
MRREVIEWER = {Juan\ M.\ Melchor},
       DOI = {10.1007/s00283-018-9788-4},
       URL = {https://doi.org/10.1007/s00283-018-9788-4},
}

@article {calkin,
    AUTHOR = {Calkin, Neil and Wilf, Herbert S.},
     TITLE = {Recounting the rationals},
   JOURNAL = {Amer. Math. Monthly},
  FJOURNAL = {American Mathematical Monthly},
    VOLUME = {107},
      YEAR = {2000},
    NUMBER = {4},
     PAGES = {360--363},
      ISSN = {0002-9890,1930-0972},
   MRCLASS = {11B75},
  MRNUMBER = {1763062},
MRREVIEWER = {H.\ L.\ Abbott},
       DOI = {10.2307/2589182},
       URL = {https://doi.org/10.2307/2589182},
}

@unpublished{wilf,
  AUTHOR = {Neil Calkin and Herbert S. Wilf},
   TITLE = {Binary partitions of integers and stern-brocot-like trees (unpublished)},
    YEAR = {1998},
     URL = {https://www.math.clemson.edu/~calkin/Papers/calkin_wilf_binary_partitions_unpublished.pdf},
}

@article {gobler,
    AUTHOR = {Gobler, Ben},
     TITLE = {Listing the rationals using continued fractions},
   JOURNAL = {Pi Mu Epsilon J.},
  FJOURNAL = {Pi Mu Epsilon Journal},
    VOLUME = {15},
      YEAR = {2022},
    NUMBER = {6},
     PAGES = {347--354},
      ISSN = {0031-952X},
   MRCLASS = {11A55},
  MRNUMBER = {4474218},
}

@article {kauffman,
    AUTHOR = {Kauffman, Louis H. and Lambropoulou, Sofia},
     TITLE = {On the classification of rational tangles},
   JOURNAL = {Adv. in Appl. Math.},
  FJOURNAL = {Advances in Applied Mathematics},
    VOLUME = {33},
      YEAR = {2004},
    NUMBER = {2},
     PAGES = {199--237},
      ISSN = {0196-8858,1090-2074},
   MRCLASS = {57M25},
  MRNUMBER = {2074397},
MRREVIEWER = {Eriko\ Hironaka},
       DOI = {10.1016/j.aam.2003.06.002},
       URL = {https://doi.org/10.1016/j.aam.2003.06.002},
}

@article {goldman,
    AUTHOR = {Goldman, Jay R. and Kauffman, Louis H.},
     TITLE = {Rational tangles},
   JOURNAL = {Adv. in Appl. Math.},
  FJOURNAL = {Advances in Applied Mathematics},
    VOLUME = {18},
      YEAR = {1997},
    NUMBER = {3},
     PAGES = {300--332},
      ISSN = {0196-8858,1090-2074},
   MRCLASS = {57M25},
  MRNUMBER = {1436484},
       DOI = {10.1006/aama.1996.0511},
       URL = {https://doi.org/10.1006/aama.1996.0511},
}

@incollection {conway,
    AUTHOR = {Conway, J. H.},
     TITLE = {An enumeration of knots and links, and some of their algebraic
              properties},
 BOOKTITLE = {Computational {P}roblems in {A}bstract {A}lgebra ({P}roc.
              {C}onf., {O}xford, 1967)},
     PAGES = {329--358},
 PUBLISHER = {Pergamon, Oxford-New York-Toronto, Ont.},
      YEAR = {1970},
   MRCLASS = {55.20},
  MRNUMBER = {258014},
MRREVIEWER = {H.\ E.\ Debrunner},
       URL = {https://doi.org/10.1016/b978-0-08-012975-4.50034-5},
}

@misc{margalit,
	AUTHOR = {Margalit, Dan},
	 TITLE = {Juggling numbers and geometry},
      YEAR = {2022},
      NOTE = {WPI {L}evi {L}. {C}onant {L}ecture, \url{https://echo360.org/media/2ec3b090-5d3e-4be9-a637-7f56b202a554/public}},
}

@misc{euclid,
     AUTHOR = {Heiberg, J. L.},
     EDITOR = {Fitzpatrick, Richard},
 TRANSLATOR = {Fitzpatrick, Richard},
      TITLE = {Euclid's elements of geometry},
       YEAR = {2008},
        URL = {https://farside.ph.utexas.edu/Books/Euclid/Elements.pdf},
}

@book {adams,
    AUTHOR = {Adams, Colin C.},
     TITLE = {The knot book: an elementary introduction to the mathematical theory of knots},
 PUBLISHER = {Providence (RI): American Mathematical Society},
      YEAR = {2004},
     PAGES = {xiv+307},
      ISBN = {0-8218-3678-1},
   MRCLASS = {57M25 (57M15 57M50 57Q45)},
  MRNUMBER = {2079925},
}

@article {brooks,
    AUTHOR = {Brooks, R. L. and Smith, C. A. B. and Stone, A. H. and Tutte,
              W. T.},
     TITLE = {The dissection of rectangles into squares},
   JOURNAL = {Duke Math. J.},
  FJOURNAL = {Duke Mathematical Journal},
    VOLUME = {7},
      YEAR = {1940},
     PAGES = {312--340},
      ISSN = {0012-7094,1547-7398},
       DOI = {10.1215/S0012-7094-40-00718-9},
   MRCLASS = {48.0X},
  MRNUMBER = {3040},
MRREVIEWER = {P.\ Scherk},
       URL = {https://doi.org/10.1007/978-0-8176-4842-8_6},
}

@article{stewart,
    AUTHOR = {Stewart, Ian},
     TITLE = {Squaring the square},
   JOURNAL = {Sci. Amer.},
  FJOURNAL = {Scientific American},
    VOLUME = {277},
      YEAR = {1997},
    NUMBER = {1},
     PAGES = {94--96},
      ISSN = {00368733, 19467087},
       DOI = {10.1038/scientificamerican0797-94},
       URL = {http://www.jstor.org/stable/24995837},
}

@misc{history,
       AUTHOR = {Anderson, Stuart E.},
        TITLE = {Tiling by squares: history and theory},
         YEAR = {2012},
         NOTE = {\url{http://www.squaring.net/history_theory/history_theory.html}}
}

\par\vspace{2\baselineskip}\noindent\textit{School of Mathematical and Statistical Sciences, Clemson University, Clemson, SC 29634, USA}

\noindent \textit{Email address:} \texttt{\{calkin, egallag, bgobler\}@clemson.edu}

\end{document}